\documentclass[11pt]{article}
\usepackage{amsfonts, latexsym, amssymb,amsmath,amstext}
\usepackage[dvips]{graphics,color}
\usepackage[mathscr]{eucal}
\newenvironment{proof}{\par\noindent{\bf Proof.}\ }{\hfill$\Box$\par\medskip}
\newtheorem{theorem}{Theorem}[section]
\newtheorem{lemma}[theorem]{Lemma}

\newcommand{\equal}{&\!\!\!=\!\!\!&}

\newcommand{\CC}{\mathbb C}
\newcommand{\PP}{\mathbb P}
\newcommand{\RR}{\mathbb R}
\newcommand{\QQ}{\mathbb Q}
\newcommand{\ZZ}{\mathbb Z}

\makeatletter
\@addtoreset{equation}{section}

\makeatother

\title{Computation of $RS$-pullback transformations for algebraic Painlev\'e VI solutions}
\author{Raimundas Vidunas\thanks{Supported by the 21 Century COE Programme
"Development of Dynamic Mathematics with High Functionality" of the Ministry
of Education, Culture, Sports, Science and Technology of Japan. 
E-mail: {\sf rvidunas@gmail.com}}\; 
and\; Alexander~V.~Kitaev\thanks{Supported by
JSPS grant-in-aide no.~$14204012$. 
E-mail: {\sf kitaev@pdmi.ras.ru}}\\
Department of Mathematics, Kyushu University, 812-8581 Fukuoka, Japan\footnotemark[1]\\
Steklov Mathematical Institute, Fontanka 27, St. Petersburg, 191023, Russia\footnotemark[2]\\
and\\
School of Mathematics and Statistics, University of Sydney,\\
Sydney, NSW 2006, Australia\footnotemark[1]\;\,\footnotemark[2]}

\date{} 

\begin{document}

\maketitle

\begin{abstract}
Algebraic solutions of the sixth Painlev\'e equation can be computed using pullback transformations 
of hypergeometric equations with respect to specially ramified rational coverings. 
In particular, as was noticed by the second author and Doran, some algebraic solutions can be constructed from a rational covering alone, without computation of the pullbacked Fuchsian equation.
But the same covering can be used to pullback different hypergeometric equations,
resulting in different algebraic Painlev\'e VI solutions.
This paper presents  
computations of explicit $RS$-pullback transformations, 
and derivation of algebraic Painlev\'e VI solutions from them. 
As an example, we present computation of all seed solutions for pull-backs of hyperbolic 
hypergeometric systems.
\vspace{24pt}\\
{\bf 2000 Mathematics Subject Classification}: 34M55, 33E17, 57M12. 
\vspace{24pt}\\
{\bf Short title}: {$RS$-pullback transformations}\\
{\bf Key words}:  $RS$-pullback transformation, isomonodromic Fucshian system, the sixth Painlev\'e equation,  algebraic solution.
\end{abstract}

\newpage

\section{Introduction}

The sixth Painlev\'e equation is, canonically,
\begin{eqnarray}
 \label{eq:P6}
\frac{d^2y}{dt^2}&=&\frac 12\left(\frac 1y+\frac 1{y-1}+\frac 1{y-t}\right)
\left(\frac{dy}{dt}\right)^2-\left(\frac 1t+\frac 1{t-1}+\frac 1{y-t}\right)
\frac{dy}{dt}\nonumber\\
&+&\frac{y(y-1)(y-t)}{t^2(t-1)^2}\left(\alpha+\beta\frac t{y^2}+
\gamma\frac{t-1}{(y-1)^2}+\delta\frac{t(t-1)}{(y-t)^2}\right),
\end{eqnarray}
where $\alpha,\,\beta,\,\gamma,\,\delta\in\CC$ are parameters. As well-known
\cite{JM}, its solutions
define isomonodromic deformations (with respect to $t$) of the $2\times2$
matrix Fuchsian equation with 4 singular points ($\lambda=0,1,t$, and
$\infty$):
\begin{equation}
 \label{eq:JM}
\frac{d}{dz}\Psi=\left(\frac{A_0}z+
\frac{A_1}{z-1}+\frac{A_t}{z-t}\right)\Psi,\qquad
\frac{d}{dz}A_k=0\quad\mbox{for } k\in\{0,1,t\}.
\end{equation}
The standard correspondence is due to Jimbo and Miwa \cite{JM}. We choose
the traceless normalization of (\ref{eq:JM}), so we assume that the
eigenvalues of $A_0$, $A_1$, $A_t$ are, respectively, $\pm\theta_0/2$,
$\pm\theta_1/2$, $\pm\theta_t/2$, and that the matrix
$A_\infty:=-A_1-A_2-A_3$ is diagonal with the diagonal entries
$\pm\theta_\infty/2$.  Then the corresponding Painlev\'e equation has the
parameters
\begin{equation}
 \label{eq:para}
\alpha=\frac{(\theta_\infty-1)^2}2,\quad
\beta=-\frac{{\theta}_0^2}2,\quad\gamma=\frac{{\theta}_1^2}2,
\quad\delta=\frac{1-{\theta}_t^2}2.
\end{equation}
We refer to the numbers $\theta_0$, $\theta_1$, $\theta_t$ and
$\theta_\infty$ as {\em local monodromy differences}.

For any numbers $\nu_1,\nu_2,\nu_t,\nu_\infty$, let us denote by
$P_{VI}(\nu_0,\nu_1,\nu_t,\nu_\infty;t)$ the Painlev\'e VI equation for the
local monodromy differences $\theta_i=\nu_i$ for $i\in\{0,1,t,\infty\}$, via
(\ref{eq:para}). Note that changing the sign of $\nu_0,\nu_1,\nu_t$ or
$1-\nu_\infty$ does not change the Painlev\'e equation. Fractional-linear
transformations for the Painlev\'e VI equation permute the 4 singular points
and the numbers $\nu_0,\nu_1,\nu_t,1-\nu_\infty$.

Similarly, for any numbers $\nu_1,\nu_2,\nu_t,\nu_\infty$ and a solution $y(t)$ of
$P_{VI}(\nu_0,\nu_1,\nu_t,\nu_\infty;t)$, let us denote by $E(\nu_0,\nu_1,\nu_t,\nu_\infty;y(t);z)$
a Fuchsian equation (\ref{eq:JM}) corresponding to $y(t)$ by the Jimbo-Miwa correspondence.
The Fuchsian equation is determined uniquely up to conjugation of $A_0,A_1,A_t$ by
a diagonal matrix (dependent  on $t$ only). In particular, \mbox{$y(t)=t$} can be considered as a
solution of $P_{VI}(e_0,e_1,0,e_\infty;t)$. The equation $E(e_0,e_1,0,e_\infty;t;z)$ is a Fuchsian equation with 3 singular points, actually without the parameter $t$. Its solutions can be expressed
in terms of Gauss hypergeometric function see \cite{JM} or the Appendix in \cite{KV3}.
We refer to $E(e_0,e_1,0,e_\infty;t;z)$ as a {\em matrix hypergeometric equation}, and
see it as a matrix form of Euler's ordinary hypergeometric equation.

We consider pullback transformations of $2\times 2$ Fuchsian systems $d\Psi(z)/dz=M(z)\Psi(z)$. 
They have the following general form:
\begin{equation}  \label{eq:rstrans}
z\mapsto R(x),\qquad \Psi(z)\mapsto S(x)\,\Psi(R(x)),
\end{equation}
where $R(x)$ is a rational function of $x$, and $S(x)$ is a Schlesinger transformation, usually designed to remove apparent singularities. For transformations to parametric 
isomonodromic equations, $R(x)$ and $S(x)$ may depend algebraically on parameter(s) as well. 
The transformed equation is
\begin{equation} \label{eq:pbacked}
\frac{d\Psi(x)}{dx}=\left(\frac{dR(x)}{dx}\,S^{-1}(x)M(R(x))S(x)-S^{-1}(x)\frac{dS(x)}{dx}\right) \Psi(x).
\end{equation}
In \cite{K2}, \cite{HGBAA}, these pullback transformations are called {\em $RS$-transformations}, meaning that they are compositions of a rational change of the independent variable $z\mapsto R(x)$
and the Schlesinger transformation $S(x)$. The Schlesinger transformation $S(x)$ is analogous here to
{\em projective equivalence} transformations $y(x)\to\theta(x)y(x)$ of ordinary differential equations.
To merge terminology, we refer to these pullback transformations as {\em $RS$-pullbacks},
or {\em $RS$-pullback transformations}.
If $S(x)$ is the identity transformation, we have a {\em direct pullback} of a Fuchsian equation.

The subject of this article is construction of $RS$-pullback transformations of matrix hypergeometric equations to isomonodromic Fuchsian systems with 4 singular points. To have so few singular points
of the transformed equation, we usually have to start with a matrix hypergeometric equation with restricted local monodromy differences, and the $R$-part $R(x)$ must define a specially ramified 
covering of $\PP^1$. In particular, the covering usually may ramify only above the 3 singular points of the hypergeometric equation, except that there is one additional simple (i.e., order 2) ramification point is allowed. Coverings ramified over 4 points of $\PP^1$ in this way are called here 
{\em almost Belyi coverings}. Recall that a {\em Belyi function} is a rational function on an algebraic
curve with at most 3 critical values; the respective covering of $\PP^1$
by the algebraic curve is ramified above a set of 3 points only.

Suitable starting hypergeometric equations and ramification patterns of almost Belyi coverings can be classified rather easily \cite{HGBAA}, \cite{Do}. This is similar to classification of algebraic transformations of Gauss hypergeometric functions \cite{V}, \cite{hhpgtr}, where Belyi functions typically occur. The computationally hard problem is construction of almost Belyi coverings from a priori suitable ramification patterns.
This leads us towards Grothendieck's theory of dessins d'enfant. In particular, Hurwitz spaces for
almost Belyi coverings with a fixed ramification pattern define isomonodromy parameters 
for the pullbacked Fuchsian equations. Effective computations
of high degree almost Belyi coverings are presented in \cite{KV2}.
In this paper, we use three coverings computed in \cite{HGBAA}.  

Computation of $S$-parts of suitable $RS$-transformations does not look hard in principle.
However, this problem is not as straightforward as finding suitable projective equivalence
transformations for scalar differential equations. General Schlesinger transformations can be constructed by composing several simple Schlesinger transformations (each shifting just two local monodromy differences), as was done in \cite{AK1}, \cite{AK2}, \cite{JM}. 
More effectively, the method in \cite{KV3} constructs Schlesinger transformations in one go,
avoiding factorization of high degree polynomials when shifting local monodromy differences
at all conjugate roots by the same integer.  In the context of isomonodromy problems,
this approach is adopted in \cite{FN} as well.

$RS$-pullback transformations to isomonodromic Fuchsian systems with 4 singular points gives 
solutions of the sixth Painlev\'e equations that are algebraic, because those solutions are determined algebraically by matrix entries of pullbacked equations (\ref{eq:pbacked}) while
those entries are algebraic functions in $x$ and the isomonodromy parameter. 
The second author conjectured in \cite{HGBAA} that all algebraic solutions of the sixth Painlev\'e equation can be obtained by $RS$-pullback transformations of matrix hypergeometric equations, up to Okamoto transformations \cite{O}. This conjecture is certainly true if the monodromy group of the
Fuchsian systems is finite, due to celebrated Klein's theorem \cite{Klein84}.
Richard Fuchs \cite{F2} soon considered extension of Klein's theorem to algebraic solutions of Painlev\'e equations. Recently, Ohyama  and Okumura \cite{Oh} showed that algebraic solutions of Painlev\'e equations from the first to the fifth do arise from pull-back transformations of confluent hypergeometric equations, affirming the formulation of R.~Fuchs.

The pullback method for computation of algebraic Painlev\'e VI solutions
was previously suggested in \cite{HGBAA}, \cite{AK2},  \cite{K2}, \cite{Do}.
This method is substantially different from the representation-theoretic approach 
of Dubrovin-Mazzocco \cite{DM} and Boalch \cite{Bo1}, \cite{Bo2}. 
Recently, \cite{LosTyk} used the representation-theoretic method 
to complete classification of algebraic Painlev\'e VI solutions.
The mentioned conjecture in \cite{HGBAA} is still interesting as a generalization of Klein's theorem.
There is a similar situation with classification of algebraic solutions of the Lam\'e equation, where representation-theoretic methods (as in \cite{BeukersWaall}) compete with Klein's pullback method
(as in \cite{Litcanu}).

One important observation is that the same rational covering $R(x)$ can be used in
several $RS$-pullback transformations. For example, here we apply the same degree 10 covering to
pullback three different matrix hypergeometric equations $E(1/7,1/2,0,1/3;t;z)$, $E(2/7,1/2,0,1/3;t;z)$
and $E(3/7,1/2,0,1/3;t;z)$. We obtain Painlev\'e solutions of, res\-pectively, $P_{VI}(1/7,1/7,1/7,2/3;t)$,
$P_{VI}(2/7,2/7,2/7,1/3;t)$ and $P_{VI}(3/7,3/7,3/7,2/3;t)$, unrelated by fractional-linear or Okamoto transformations. The first Painlev\'e solution is a fractional-linear version of solution 
\cite[(3.16)--(3.17)]{K4}.  The second Painlev\'e solution is the same as in \cite[page 106]{Bo3}.
The third Painlev\'e solution is new.

The article is organized as follows. Section \ref{sec:pbmaps} presents the covering of degree 10
for our exmaples; it was previously used in \cite{K4}. There we also mention
how some Painlev\'e VI solutions can be computed from the rational coverings alone, without computation of full $RS$-transformations. This kind of possibility is noticed 
in  \cite{HGBAA}, \cite{Do}, 
and is summarized in Theorem \ref{kit:method} below. 
In Section \ref{sec:fullrs2}  a more general Theorem \ref{th:llentry} from \cite{KV3} is cited.
Thereby a direct formula for algebraic Painlev\'e VI solutions is given, 
with minimum information from full $RS$-transformations. 
In Section \ref{sec:algsols}, representative $RS$-pullback transformations 
of ``hyperbolic" hypergeometric equations 
$E(1/2,1/3,0,1/7;t;z)$ and $E(1/2,1/3,0,1/8;t;z)$ to isomonodromic Fuchsian systems
with 4 singularities are sumarized, and the corresponding Painlev\'e VI solutions are presented
(hereby complementing \cite{K4}).  The Appendix 
presents a formula for composition of two quadratic transformations
of Painlev\'e VI solutions; a general degree formula 
for the almost Belyi coverings relevant to algebraic Painlev\'e VI solutions;
and geometric interpretation of the latter formula.

The authors prepared {\sf Maple 9.5} worksheets supplementing this article and \cite{KV2},  \cite{KV3},
with the formulas in {\sf Maple} input format, and demonstration of key
computations. To access the worksheet, readers may contact the authors,
or search a current website of the first author on the internet.

\section{The working covering and $RS$-transformations}
\label{sec:pbmaps}

First we introduce notation for ramification patterns, and later for $RS$-transformations.
A ramification pattern for an almost Belyi covering of degree $n$ is denoted by
$R_4(P_1|\,P_2|\,P_3)$, where $P_1,P_2,P_3$ are three partitions of $n$
specifying the ramification orders above three points. The ramification pattern
above the fourth ramification locus is assumed to be $2+1+1+\ldots+1$.
By {\em the extra ramification point} we refer to the simple ramification point in the fourth fiber.
The Hurwitz space for such a ramification pattern is generally one-dimensional 
\cite[Proposition 3.1]{zvonkin}. 

We use only genus 0 almost Belyi coverings, and write them as $\PP_x^1\to\PP_z^1$, 
meaning that the projective line with the projective coordinate $x$ is mapped to the projective line
with the coordinate $z$. Then the total number
of parts in $P_1$, $P_2$, $P_3$ must be equal to $n+3$, according to \cite[Proposition 2.1]{HGBAA};
this is a consequence of Riemann-Hurwitz formula.

The similar notation for a ramification
pattern for a Belyi function is $R_3(P_1|\,P_2|\,P_3)$, as in \cite{AK1}, \cite{K4}.
The total number of parts in $P_1$, $P_2$, $P_3$ must be equal to $n+2$,
as stated in \cite[Lemma 2.4]{V} or \cite[Proposition 2.1]{HGBAA}.

Our working almost Belyi covering has the following ramification pattern:
\begin{eqnarray} \label{eq:ramif7}
& R_4\big(7+1+1+1\;|\;2+2+2+2+2\;|\;3+3+3+1\big).
\end{eqnarray}
The covering has degree 10.  
The three specified fibers with ramified points
can be brought to any three distinct locations 
by a fractional-linear transformation of $\PP^1_z$. We assign the first partition to $z=0$,
and the next two partitions --- to $z=1$ and $z=\infty$ respectively.
Similarly, by a fractional-linear transformation of $\PP^1_x$ we may choose any three 
$x$-points\footnote{Strictly speaking, the $x$-points in our settings are curves, or branches, parametrized by an isomonodromy parameter $t$ or other parameter, since the Hurwitz spaces
for almost Belyi maps are one-dimensional. For simplicity, we ignore the dimensions introduced
by  such parameters, and consider a one-dimensional Hurwitz space as a generic point.}
as $x=0$, $x=1$, $x=\infty$. 

All coverings with ramification pattern (\ref{eq:ramif7}) can be computed 
on modern computers either using the most straightforward method, or an improved method \cite{KV2}
that uses differentiation. 
Up to fractional-linear transformations and reparametrization,
there is one general such covering given by
\begin{equation} \label{eq:phi10}
\varphi_{10}(x)=\frac{x^7\,F_{10}}{4\,G_{10}^3}, \qquad\mbox{or}\qquad
\varphi_{10}(x)-1=\frac{P_{10}^2}{4\,G_{10}^3},
\end{equation}
where 
\begin{eqnarray} \label{eq:ph10p}
F_{10}\equal 9s^2x^3-2(2s^3+6s^2+15s-16)x^2+3(8s^2+8s-13)x-36(s-1), \nonumber\\
G_{10}\equal 2(s+1)x^3-(s^2+4s+10)x^2+6(s+2)x-9,\\  \nonumber
P_{10}\equal 3sx^5-3(2s^2+6s+7)x^4+2(s^3+6s^2+30s+35)x^3\\ \nonumber
&&-18(s^2+4s+7)x^2+54(s+2)x-54.
\end{eqnarray}
The extra ramification point is $x=7(s-1)/s(s+1)$. 

For direct applications to the Painlev\'e VI equation, it is required to normalize
the point above $z=\infty$ with the deviating ramification order 1
and the three nonramified points above $\{0,1,\infty\}\subset\PP^1_z$ as
$x=0$, $x=1$, $x=\infty$, $x=t$. We refer to explicit almost Belyi coverings normalized
this way as {\em properly normalized}. 
A properly normalized covering with ramification pattern (\ref{eq:ramif7}) was first
computed in \cite{HGBAA}. 
To get a properly normalized expression, we reparametrize 
\begin{equation} \label{eq:ts10}
s=-\frac{(u+2)(u^2-u+2)}{2(u-1)},
\end{equation}
and make the fractional-linear transformation
\begin{equation} \label{eq:lambd10}
x\mapsto \frac{(u-1)^2(u^2+3)w}{9(u^2-u+2)^2}(2x-1)
-\frac{(u-1)(u^5+u^4-2u^3+18u^2-9u+27)}{9(u^2-u+2)^2},
\end{equation}
where $w=\sqrt{(u-1)(u+5)(u^2+3)}$. The obtained expression is
\begin{eqnarray}
\hspace{-18pt} \widehat{\varphi}_{10}(x)=
-\frac{(u-1)^2(u+2)^2w^3}{8(u\!+\!5)(u^3\!+\!u^2\!-\!2u\!+\!6)^3}
\frac{x\,(x-1)\left(x-t_{10}\right)\left(x-t^*_{10}\right)^7 }
{\left( \left(x^2-x\right)(x-\frac12-L_3) -L_4(x-\frac12)+L_5 \right)^3},
\end{eqnarray}
where
\begin{eqnarray} \label{eq:t10}
t_{10}\equal \frac12+\frac{u^9+3u^8-3u^7+7u^6-21u^5+21u^4-161u^3-27u^2-144u-108}
{2(u-1)^3(u+2)^2(u^2+3)\sqrt{(u-1)(u+5)(u^2+3)}},\quad\\
t^*_{10}\equal \frac12+\frac{u^5+u^4-2u^3+18u^2-9u+27}{2(u-1)(u^2+3)\sqrt{(u-1)(u+5)(u^2+3)}},
\end{eqnarray}
and
\begin{eqnarray*}
L_3\equal \textstyle \frac{(u^5+4u^4+u^3+18u^2+24u+36)(u^7+14u^4-21u^3+42u^2+36)}
{8(u-1)^2(u+5)(u^2+3)^3(u^3+u^2-2u+6)},\\
L_4\equal \textstyle \frac{3(u^{10}-6u^8+28u^7-99u^6+252u^5-668u^4+1008u^3-1212u^2+672u-408)}
{8w(u-1)^3(u^2+3)(u^3+u^2-2u+6)},\\
L_5\equal \textstyle \!\frac{u^{15}+5u^{14}+28u^{12}+98u^{11}-126u^{10}+616u^9-184u^8+333u^7+1785u^6-1512u^5+3276u^4+6048u^2+3888u+1296}{16w(u-1)^3(u+5)(u^2+3)^3(u^3+u^2-2u+6)}.\end{eqnarray*}
The Hurwitz space parametrising this
properly normalized almost Belyi covering has still genus 0. 
To get the rational covering $\lambda_1(\lambda)$ in \cite{K4}, one has to consider 
$t_{10}\big/\widehat{\varphi}_{10}(x)$, and substitute $x\mapsto t_{10}/x$, $u\mapsto 2/s-1$.

In \cite{HGBAA}, the following symbol is introduced to denote 
$RS$-pullback transformations of $E(e_0,e_1,0,e_\infty;t;z)$ with respect to a covering with ramification
pattern $R_4(P_0|P_1|P_\infty)$: 
\begin{equation}
RS^2_4\left( \,e_0\, \atop \,P_0\, \right| 
{\,e_1\, \atop \,P_1\, } \left| \,e_\infty\, \atop \,P_\infty\, \right),
\end{equation}
where the subscripts 2 and 4 indicate a second order Fuchsian system with 4 singular points 
after the pullback. 
We assume the same assignment of the fibers $z=0$, $z=1$, $z=\infty$ as for the $R_4$-notation.
Location of the $x$-branches $0,1,t,\infty$ does not have to be normalized.
As was noticed in \cite{HGBAA} and \cite{Do}, some
algebraic Painlev\'e VI solutions determined by $RS$-pullback transformations 
$RS^2_4\left( \,1/k_0\, \atop P_0 \right| {\,1/k_1\, \atop P_1 } \left| \,1/k_\infty\, \atop P_\infty \right)$,
with $k_0,k_1,k_\infty\in\ZZ$, can be calculated from the rational covering alone, without actual
computation of the full $RS$-pullbacks. 
We discuss this possibility in Section \ref{sec:painleve}. 
Our covering $\widehat{\varphi}_{10}(x)$ immediately gives a solution of 
$P_{VI}(1/7,1/7,1/7,2/3;t)$.
In Section \ref{sec:fullrs2} we formulate a direct way to obtain algebraic Painlev\'e VI solutions
via computation of suitable syzygies between $x^2$ (or $x^3$), $P_{10}$, $G_{10}$. 
We obtain algebraic solutions of $P_{VI}(2/7,2/7,2/7,1/3;t)$ and $P_{VI}(3/7,3/7,3/7,2/3;t)$
by implicitly using $RS$-pullback transformations 
$RS^2_4\left( 2/7 \atop 7+1+1+1 \right| {1/2 \atop 2+2+2+2+2} \left| 1/3 \atop 3+3+3+1 \right)$ and
$RS^2_4\left( 3/7 \atop 7+1+1+1 \right| {1/2 \atop 2+2+2+2+2} \left| 1/3 \atop 3+3+3+1 \right)$, respectively. 

\section{Pullback coverings and algebraic Painlev\'e VI solutions}
\label{sec:painleve}


As noticed in  \cite{HGBAA} and \cite{Do}, some algebraic Painlev\'e VI solutions can be computed knowing just a pullback covering, without computation of pullbacked Fuchsian equations of full
$RS$-transformations removing all apparent singularities of a direct pullback. Here we formulate the most interesting general situation. 
\begin{theorem} \label{kit:method}
Let $k_0,k_1,k_\infty$ denote three integers, all $\ge 2$. Let $\varphi:\PP^1_x\to\PP^1_z$ denote
an almost Belyi map, dependent on a parameter $t$. Suppose that the following conditions are satisfied:
\begin{itemize}
\item[(i)] The covering $z=\varphi(x)$ is ramified above the points $z=0$, $z=1$, $z=\infty$; there is
one simply ramified point $x=y$ above $\PP_z^1\setminus\{0,1,\infty\}$; 
and there are no other ramified points.
\item[(ii)] The points $x=0$, $x=1$, $x=\infty$, $x=t$ 
lie above the set $\{0,1,\infty\}\subset\PP^1_z$.
 \item[(iii)] The points in $\varphi^{-1}(0)\setminus\{0,1,t,\infty\}$ are all ramified with the order $k_0$.
The points in $\varphi^{-1}(1)\setminus\{0,1,t,\infty\}$ are all ramified with the order $k_1$.
The points in $\varphi^{-1}(\infty)\setminus\{0,1,t,\infty\}$ are all ramified with the order $k_\infty$.
\end{itemize}
Let $a_0,a_1,a_t,a_{\infty}$ denote the ramification orders at $x=0,1,t,\infty$, respectively. 
Then the point $x=y$, as a function of $x=t$, is an algebraic solution of 
\begin{equation} \label{eq:p6kit}
P_{VI}\left(\frac{a_0}{k_{\varphi(0)}},\frac{a_1}{k_{\varphi(1)}},
\frac{a_t}{k_{\varphi(t)}}, 1-\frac{a_{\infty}}{k_{\varphi(\infty)}};t\right).
\end{equation}
\end{theorem}
\begin{proof} Let $R_4(P_0|P_1|P_\infty)$ denote the ramification pattern of the covering $z=\varphi(x)$. 
We aim for an $RS$-pullback transformation 
$RS^2_4\left( \,1/k_0\, \atop P_0 \right| {\,1/k_1\, \atop P_1 } \left| \,1+1/k_\infty\!\atop P_\infty \right)$
with respect to $\varphi(x)$. Let $d$ denote the degree of $\varphi(x)$.
For time being, we assume that the point $x=\infty$ lies above $z=\infty$.

The direct pullback of the hypergeometric equation $E(1/k_0,1/k_1,0,1+1/k_\infty;t;z)$ with respect to
$\varphi(x)$ has apparent singularities at the points mentioned in part {\em(iii)} above.
Nonapparent singularities are possibly $x=0$, $x=1$, $x=t$ and $x=\infty$. 
The lower-left entry of the direct pullback is equal, up to a factor independent of $x$,
to  $\varphi'/\varphi(1-\varphi)$, which is the logarithmic derivative of $\varphi/(\varphi-1)$.
The poles of this rational function are simple, and they are precisely the points
above $z=0$ and $z=1$. The zeroes of the rational function are the following:
the extra ramification point of $\varphi$ (a simple zero); and the points above $z=\infty$,
with multiplicities one less than the respective ramification orders.

Notice that if we apply
a Schlesinger transformation of the upper triangular form $S=\frac1{\sqrt{(x-\alpha_1)(x-\alpha_2)}}
{ x-\alpha_1 \quad\alpha_3\ \ \choose \ \ 0 \quad x-\alpha_2}$, where $\alpha_1$, $\alpha_2$, 
$\alpha_3$ are independent of $x$, 
then the lower-left entry of the matrix differential equation
changes by the factor $(x-\alpha_1)/(x-\alpha_2)$ and a factor independent of $x$.
If the point $x=\alpha_2$ is above $z=\infty$,
this Schlesinger transformation (with appropriate $\alpha_3$)  
decreases  the local monodromy differences at $x=\alpha_1$ and $x=\alpha_2$ by 1.
Similarly, the Schlesinger transformation
$S=\frac1{\sqrt{x-\alpha_1}}{ x-\alpha_1 \quad 0\ \choose \ 0 \qquad 1}$ 
changes the local monodromy differences at $x=\alpha_1$ and $x=\infty$ by 1,
and it multiplies the lower-left entry by the factor $x-\alpha_1$  (and a factor independent of $x$).

Let $h$ denote the number of distinct apparent singularities above $z=\infty$.
There are in total $(d+3)-4-h$ apparent singularities above $z=0$ and $z=1$. 
We can construct $d-1-h$ simple Schlesinger transformations of the forms presented just above,
so that $\alpha_1$ runs through the set of apparent singularities above $z=0$ and $z=1$, and
each point $x=\alpha_2$ or $x=\infty$ above $z=\infty$ is chosen $n_x$ times, where
\[ 
n_x=\left\{ \begin{array}{rl}
\mbox{the ramification order at $x$, minus 1}, &  \mbox{if $x=\infty$ or an apparent singularity};\\
\mbox{the ramification order at $x$}, &  \mbox{otherwise}. \end{array} \right.
\] 
The composite effect of these $d-1-h$ transformations is removal of all apparent singularities
above $z=0$, $z=1$, $z=\infty$; and reducing the local monodromy difference at $x=\infty$ from
$a_{\infty}+a_{\infty}/k_{\infty}$ to $1+a_{\infty}/k_{\infty}$.
The local monodromy differences at the other singularities are $a_0/k_{\varphi(0)}$,
$a_1/k_{\varphi(1)}$, $a_t/k_{\varphi(t)}$ after the composite transformation.
Hence the transformed equation has (at most) four singularities. The transformed equation is
\begin{equation} \label{kit:fuchs}
E\left(\frac{a_0}{k_{\varphi(0)}},\frac{a_1}{k_{\varphi(1)}},
\frac{a_t}{k_{\varphi(t)}}, 1+\frac{a_{\infty}}{k_{\infty}};\widetilde{y}(t);x\right),
\end{equation}
where the Painlev\'e VI solution $\widetilde{y}(t)$ is determined by lower-left entry of the transformed equation. The lower-left entry is changed from $\varphi'/\varphi(1-\varphi)$ to a rational function
whose numerator has only one root. The single root must be the extra ramification point of $\varphi(x)$. Hence $\widetilde{y}(t)$ can be identified with the branch $x=y$. 
It is a solution of  $P_{VI}\!\left(a_0/k_{\varphi(0)},{a_1}/{k_{\varphi(1)}},
{a_t}/{k_{\varphi(t)}}, 1+{a_{\infty}}/{k_{\infty}};t\right)$
which is the same equation as  (\ref{eq:p6kit}). 

If the point $x=\infty$ does not lie above $z=\infty$, we can move the point $z=\infty$ by the
fractional-linear transformations. That would only permute the three fibers, and change the rational 
function $\varphi$ to $1/\varphi$, $1/(1-\varphi)$, $1-1/\varphi$ or $\varphi/(\varphi-1)$ . 
Action of fractional-linear transformations on local monodromy differences is compatible
with the form (\ref{eq:p6kit}).
\end{proof}

The above theorem is  a special case of \cite[Theorem 2.1]{HGBAA}, with all $k_j^i$'s equal to $1$,
and with correct parameters in (\ref{eq:p6kit}).  
Theorem 4.5 in \cite{Do} is a more general statement, but without identification of transformed local monodromy differences.

In \cite{HGBAA}, it is regularly implied that the Painlev\'e VI solutions obtained with 
Theorem \ref{kit:method} arise from $RS$-pullback transformations of the type
$RS^2_4\left( \,1/k_0\, \atop P_0 \right| {\,1/k_1\, \atop P_1 } \left| \,1/k_\infty\!\atop P_\infty \right)$.
However, the above proof actually uses transformation 
$RS^2_4\left( \,1/k_0\, \atop P_0 \right| {\,1/k_1\, \atop P_1 } \left| \,1+1/k_\infty\!\atop P_\infty \right)$. 
On the other hand, it is apparent from classification \cite{K4} of rational coverings for 
$RS^2_4$-pullback transformations relevant to the sixth Painlev\'e equation that either
$k_0=2$ or $k_1=2$ or $k_\infty=2$.  Once we assume $k_\infty=2$, the transformations
types 
$RS^2_4\left( \,1/k_0\, \atop P_0 \right| {\,1/k_1\, \atop P_1 } \left| \,1/k_\infty\!\atop P_\infty \right)$
and
$RS^2_4\left( \,1/k_0\, \atop P_0 \right| {\,1/k_1\, \atop P_1 } \left| \,1+1/k_\infty\!\atop P_\infty \right)$
are the same or related by extra Schlesinger transformations. If $k_0=2$ or $k_1=2$,
we still can relate the two transformation types via Schlesinger transformations.
Hence, the $RS$-pullback transformation implied in Theorem \ref{kit:method} can be realized 
as $RS^2_4\left( \,1/k_0\, \atop P_0 \right| {\,1/k_1\, \atop P_1 } \left| \,1/k_\infty\!\atop P_\infty \right)$
as well.


Application of Theorem \ref{kit:method} to $\widehat{\varphi}_{10}(x)$ gives
this  solution 
of  $P_{VI}(1/7,1/7,1/7,2/3;t_{10})$:
\begin{eqnarray} \label{eq:yt71}
y_{71}\equal\frac12+\frac{(u+5)(u^6-u^5+3u^4-13u^3+4u^2-18u-12)}
{2(u-1)(u+2)(u^3+u^2-2u+6)\sqrt{(u-1)(u+5)(u^2+3)}}.
\end{eqnarray}
A parametrization of $t_{10}$ is given in (\ref{eq:t10}).
To get the solution of $P_{VI}(1/3,1/7,1/7,6/7;t_{10})$ 
in \cite[(3.6)--(3.7)]{K4}, one has to consider the function $t_{10}/y_{71}$ and
substitute   $u\mapsto 2/s-1$. Our implied $RS$-transformation is
$RS^2_4\left( 1/7 \atop 7+1+1+1 \right|  {1/2 \atop 2+2+2+2+2} \left| 1/3 \atop 3+3+3+1 \right)$.

\section{Painlev\'e solutions from more general $RS$-pullback transformations}
\label{sec:fullrs2}

By the Jimbo-Miwa correspondence, a Painlev\'e VI solution is determined by the lower-left entry of 
a pullbacked Fuchsian system. By the results in \cite[Section 4]{KV3}, that lower-left entry is determined by a syzygy $(U_2,V_2,W_2)$ between $F$, $G$, $H$; that is, a polynomial solution of $FU_2+GV_2+HW_2=0$. If the shift $\delta$ of local monodromy differences at $x=\infty$ is small,
that syzygy is determined by degree bounds of its components.
The following theorem summarizes the situation. 
\begin{theorem} \label{th:llentry}
Let $z=\varphi(x)$ denote a rational covering, 
and let $F(x)$, $G(x)$, $H(x)$ denote polynomials in $x$. 
Let $k$ denote the order of the pole of $\varphi(x)$ at $x=\infty$.
Suppose that the direct pullback of $E(e_0,e_1,0,e_\infty;t;z)$
with respect to $\varphi(x)$ is a Fuchsian equation with the following singularities:
\begin{itemize}
\item Four singularities are $x=0$, $x=1$, $x=\infty$ and $x=t$, with the local monodromy differences 
$d_0$, $d_1$, $d_t$, $d_\infty$, respectively. The point $x=\infty$ lies above $z=\infty$.
\item All other singularities in $\PP^1_x\setminus\{0,1,t,\infty\}$ are apparent singularities.
The apparent singularities above $z=0$ (respectively, above $z=1$, $z=\infty$) are 
the roots of $F(x)=0$ (respectively, of $G(x)=0$, $H(x)=0$).
Their local monodromy differences are equal to the multiplicities of those roots. 
\end{itemize}
Let us denote $\Delta=\deg F+\deg G+\deg H$, and let $\delta\le\max(2,k)$ 
denote a non-negative integer such that $\Delta+\delta$ is even.
Suppose that $(U_2,V_2,W_2)$ is a syzygy between the three polynomials $F$, $G$, $H$, 
satisfying, if $\delta=0$,
\begin{equation} \label{eq:llsd0} \textstyle
\deg U_2=\frac{\Delta}2-\deg F,\qquad \deg V_2=\frac{\Delta}2-\deg G,\qquad
\deg W_2<\frac{\Delta}2-\deg H, 
\end{equation}
or, if $\delta>0$,
\begin{equation} \label{eq:llsd1} \textstyle
\deg U_2<\frac{\Delta+\delta}2-\deg F, \quad
\deg V_2<\frac{\Delta+\delta}2-\deg G,\quad
\deg W_2=\frac{\Delta-\delta}2-\deg H.
\end{equation}
Then the numerator of the (simplified) rational function 
\begin{eqnarray} \label{eq:llentry}
\frac{U_2W_2}{G}\!\left(\! \frac{(e_0-e_1+e_\infty)}2\frac{\varphi'}{\varphi}
-\frac{(FU_2)'}{FU_2}+\frac{(HW_2)'}{HW_2} \right) 
+\frac{(e_0-e_1-e_\infty)}2\frac{V_2W_2}{F}\frac{\varphi'}{\varphi-1} \nonumber\\
+\frac{(e_0+e_1-e_\infty)}2\frac{U_2V_2}{H}\frac{\varphi'}{\varphi\,(\varphi-1)},
\end{eqnarray}
has degree $1$ in $x$, and the $x$-root of it 
is an algebraic solution of  $P_{VI}(d_0,d_1,d_t,d_\infty+\delta;t)$.
\end{theorem}
\begin{proof} See Theorem 5.1 in \cite{KV3}.
\end{proof}

Alternative forms of expression (\ref{eq:llentry}) are given in formulas (5.17)--(5.22) in \cite{KV3}.
For greater $\delta$, formula (\ref{eq:llentry}) is still valid for a suitable syzygy $(U_2,V_2,W_2)$,
but that syzygy depends on initial coefficients of local solutions at $z=0$ of the original hypergeometric equation. Taking only small shifts $\delta<\max(2,k)$ at $x=\infty$ seems to be enough to generate interesting ``seed" solutions of the sixth Painlev\'e equation. 

We can apply this theorem to obtain algebraic solutions of $P_{VI}(1/7,1/7,1/7,2/3;t)$,
$P_{VI}(2/7,2/7,2/7,1/3;t)$ and $P_{VI}(3/7,3/7,3/7,2/3;t)$. Implicitly, we apply
pullback transformations 
$RS^2_4\left( 1/7 \atop 7+1+1+1 \right| {1/2 \atop 2+2+2+2+2} \left| 1/3 \atop 3+3+3+1 \right)$,
$RS^2_4\left( 2/7 \atop 7+1+1+1 \right| {1/2 \atop 2+2+2+2+2} \left| 1/3 \atop 3+3+3+1 \right)$ and
$RS^2_4\left( 3/7 \atop 7+1+1+1 \right| {1/2 \atop 2+2+2+2+2} \left| 1/3 \atop 3+3+3+1 \right)$, respectively.  Like in Section \ref{sec:painleve}, we work with the covering $z=\varphi_{10}(x)$ rather
than with the normalized covering $z=\widehat{\varphi}_{10}(x)$ while computing syzygies,  
and apply reparametrization (\ref{eq:ts10}) and normalizing fractional-linear transformation 
(\ref{eq:lambd10}) at the latest stage. We have $k=1$.
Therefore recall the definition of $F_{10}$, $G_{10}$ and $P_{10}$ in (\ref{eq:ph10p}).
We take $\delta=0$ for the second $RS$-transformation, or $\delta=1$ for the other two.
We have to compute syzygies between $F=x$ (or, respectively, $F=x^2$, or $F=x^3$) 
and $G=P_{10}$, $H=G_{10}$. 

The syzygy for a solution of $P_{VI}(1/7,1/7,1/7,2/3;t)$ is $(G_{10},0,-x)$, up to a scalar multiple.
With this trivial syzygy, the solution is the same $\varphi_{10}(t_{10})$ as in (\ref{eq:yt71}).
In fact, Theorem \ref{th:llentry} reduces to Theorem \ref{kit:method} whenever one of syzygy components is zero; see \cite[Remark 5.2]{KV3}.  

The full $RS$-pullback 
$RS^2_4\left( 1/7 \atop 7+1+1+1 \right| {1/2 \atop 2+2+2+2+2} \left| 1/3 \atop 3+3+3+1 \right)$ 
would give a solution $\widetilde{y}_{71}(t_{70})$ of $P_{VI}(1/7,1/7,1/7,-2/3;t_{10})$ as well.
The equation $P_{VI}(1/7,1/7,1/7,8/3;t_{10})$ is identical. 
It turns out that the same Painlev\'e solution can be obtained by applying 
Theorem \ref{th:llentry} with $\delta=3$. (Have a look at the second part of  \cite[Remark 5.3]{KV3}.)
However, since $\delta=3>\max(2,1)$ we are not given restrictions on the syzygy $(U_2,V_2,W_2)$, 
and additional knowledge of the normalized solutions of $E(1/7,1/2,0,1/3;t;z)$ at $z=\infty$ is needed. 
The syzygy can be eventually computed to be
\begin{eqnarray*}
\left( -63s^2x^4+(74s^3+222s^2+285s-52)x^3-2(8s^4+48s^3+257s^2+297s-130)x^2\;\right.\\
\left.+6(16s^3+64s^2+101s-52)x-144s^2-288s+234,\; 21s,\; 26(s+1)^2x-126s \right).
\end{eqnarray*}
The numerator of simplified expression (\ref{eq:llentry}) is then indeed linear in $x$.  
The solution $\widetilde{y}_{71}(t_{10})$ is rather stupendous:{\footnotesize
\begin{eqnarray*}
\frac12+\frac{(u+5)(65u^{18}+195u^{17}-195u^{16}+325u^{15}
-1104u^{14}+\ldots
-248931u^2-299835u+222534)}{10(u+2)\sqrt{(u-1)^3(u+5)(u^2+3)}(13u^{15}+65u^{14}+42u^{11}
-1050u^{10}+\ldots
-37611u^2+63927u-783)}.
\end{eqnarray*}}On 
the other hand, this solution can be obtained by applying a series of Okamoto transformations
to $y_{71}(t_{10})$.

To get a solution of $P_{VI}(2/7,2/7,2/7,1/3;t)$ we apply Theorem \ref{th:llentry} 
with $(F,G,H)=\left(x^2,P_{10},G_{10}\right)$. With $\delta=0$, the degree specifications 
in (\ref{eq:llsd0}) are 
\begin{equation}  \label{eq:s5degb}
\deg U_2=3,  \qquad \deg V_2=0,\qquad \deg W_2<2.
\end{equation}
As expected, there is one syzygy satisfying these bounds, up to a constant multiple:
\begin{eqnarray} \label{eq:syzs32} \textstyle 
\left( 3sx^3-(2s^2+6s+13)x^2+6(2s+3)x-18, -1, -2(s+2)x+6 \right).
\end{eqnarray}
With this syzygy, expression (\ref{eq:llentry}) is equal to
\begin{equation}
\frac{4\left(s(2s^2+4s-19)x-3(2s^2-12s+7)\right)}
{7\,F_{10}}.
\end{equation}
The form is as expected: the numerator has degree 1 in $x$, while the denominator 
is a cubic polynomial in $x$. After reparametrization (\ref{eq:ts10}) and 
normalizing fractional-linear transformation
(\ref{eq:lambd10}) the denominator polynomial surely factors as $x(x-1)(x-t_{10})$, 
with $t_{10}$ given in (\ref{eq:t10}). The $x$-root of the transformed numerator gives
the following solution $y_{72}(t_{10})$ of $P_{VI}(2/7,2/7,2/7,1/3;t_{10})$:
\begin{eqnarray} \label{eq:yt72}
\!\! y_{72}\equal\frac12+\frac{(u+5)(u^8+u^7+u^6-u^5+8u^4-82u^3-54u^2-90u-108)}
{2(u\!+\!2)(u^6\!+\!2u^5\!-\!3u^4\!+\!8u^3\!-\!26u^2\!+\!60u\!-\!6)\sqrt{(u\!-\!1)(u\!+\!5)(u^2\!+\!3)}}.
\end{eqnarray}
To relate to Boalch's parametrization in \cite[page 106]{Bo3} for the same solution, 
we have to substitute $u\mapsto (s+5)/(s-1)$ into the expressions for $y_{72}$ and $t_{10}$.

A solution $\widetilde{y}_{72}(t_{10})$ of $P_{VI}(2/7,2/7,2/7,-1/3;t_{10})$ can be computed 
without extra knowledge of the normalized solutions at $z=\infty$. 
The identical Painlev\'e equation is $P_{VI}(2/7,2/7,2/7,7/3;t_{10})$,
and Theorem \ref{th:llentry} can be applied with $\delta=2$.  
The following syzygy fits into formula (\ref{eq:llentry}): 
\begin{eqnarray*}
\left( -69s(s+1)x^3+(32s^3+128s^2+325s-65)x^2-6(32s^2+59s-15)x+288s-90, \right.\\
\left.-5s-5, 42sx^2-10(s+1)(s+2)x+30+30s\right).
\end{eqnarray*}

Application of Theorem \ref{th:llentry} with $(F,G,H)=\left(x^3,P_{10},G_{10}\right)$ and $\delta=1$
gives a solution of $P_{VI}(3/7,3/7,3/7,2/3;t)$. The degree bounds are
$\deg U_2<3$, $\deg V_2<1$, $\deg W_2=2$.
An appropriate syzygy is
\begin{equation}
\left( -(s+4)x^2+(2s+7)x-6, -1, 2x^2-2(s+2)x+6 \right)
\end{equation}
Simplified expression (\ref{eq:llentry}) has the unique $x$-root $x=-(2s-5)(4s-7)/s(10s-11)$.
After reparametrization (\ref{eq:ts10}) and  normalizing fractional-linear transformation
(\ref{eq:lambd10}) we derive the following solution $y_{73}(t_{10})$ of $P_{VI}(3/7,3/7,3/7,2/3;t_{10})$:
\begin{eqnarray}
\!\!y_{73}\equal\frac12+\frac{(u+5)(5u^7-10u^6+5u^5-20u^4+13u^3-68u^2-3u-30)}
{2(u-1)^2(u+2)(5u^3+5u^2+11u+9)\sqrt{(u\!-\!1)(u\!+\!5)(u^2\!+\!3)}}.
\end{eqnarray}
This solution cannot be obtained by Okamoto, fractional-linear and quadratic transformations from previously know solutions.

\section{Pull-backs of hyperbolic hypergeometric equations}
\label{sec:algsols}

Here we survey $RS^2_4$-pullback transformations of {\em hyperbolic} hypergeometric equations
$E(e_0,e_1,0,e_\infty;t;z)$; these are defined by the properties that $1/e_0,1/e_1,1/e_\infty$ are positive integers and $e_0+e_1+e_\infty<1$. These pullback coverings (and corresponding Okamoto orbits of algebraic Painlev\'e VI solutions)  are classified in \cite{K4} and \cite{Do}. The following ramification patterns are possible:
\begin{eqnarray} \label{eq:ramif10a}
& R_4\big(7+1+1+1\;|\;2+2+2+2+2\;|\;3+3+3+1\big),\\ 
\label{eq:ramif12a}
& R_4\big(3+3+3+3\;|\;2+2+2+2+2+2\;|\;7+2+1+1+1\big),\\
\label{eq:ramif12b}
& R_4\big(3+3+3+3\;|\;2+2+2+2+2+2\;|\;8+1+1+1+1\big),\\
\label{eq:ramif18}
& R_4\big(3\!+\!3\!+\!3\!+\!3\!+\!3\!+\!3\;|\;2\!+\!2\!+\!2\!+\!2\!+\!2\!+\!2\!+\!2\!+\!2\!+\!2\;|\;
7\!+\!7\!+\!1\!+\!1\!+\!1\!+\!1\big).
\end{eqnarray}
The coverings have degree 10, 12, 12, 18, respectively.  

The generic degree 10 covering is our $\varphi_{10}(x)$, up to reparametrization. 
We already considered the solutions $y_{71}(t_{10})$,
$y_{72}(t_{10})$, $y_{73}(t_{10})$ representing three possible Okamoto orbits.

The generic degree 12 covering with ramification (\ref{eq:ramif12a}) is 
\begin{equation}
\phi_{12}(x)=\frac{4 \left(x^4+2s(3s+1)x^3+2s(5s+2)x^2+4s^2x+s^2\right)^3}
{27s(s+1)^3x^7\left(4x^3+4s(8s+5)x^2+s(13s+1)x+4s^2\right)}.
\end{equation}
It can be normalized with the substitutions
\begin{eqnarray}
&& s\mapsto-\frac{(u+1)^2(u-1)^2}{(u^2+7)(u^2+u+2)(u^2-u+2)}, \\
&& x\mapsto \frac{(u+1)(u-1)^2}{2(u^2-u+2)^2}-\frac{u^3(u+1)(u-1)(u^2+3)\,x}{(u^2+u+2)^2(u^2-u+2)^2}.
\end{eqnarray}
A normalized expression for $1\big/{\phi}_{12}(x)$ 
is presented in \cite{HGBAA}, reparametrized with $u\mapsto 1/s$. Similarly as with $\varphi_{10}(x)$,
we can pullback  $E(1/3,1/2,0,1/7;t;z)$, $E(1/3,1/2,0,2/7;t;z)$ and $E(1/3,1/2,0,3/7;t;z)$ 
with respect to a properly normalized $\phi_{12}(x)$  and derive\footnote{The implied 
$RS$-pullback transformations are, respectively, 
$RS^2_4\left( 1/3 \atop 3+3+3+3 \right|  {1/2 \atop 2+2+2+2+2+2} \left| 1/7 \atop 7+2+1+1+1 \right)$,
$RS^2_4\left( 1/3 \atop 3+3+3+3 \right|  {1/2 \atop 2+2+2+2+2+2} \left| 2/7 \atop 7+2+1+1+1 \right)$ and
$RS^2_4\left( 1/3 \atop 3+3+3+3 \right|  {1/2 \atop 2+2+2+2+2+2} \left| 3/7 \atop 7+2+1+1+1 \right)$.  
As indicated in \cite{K4}, one may also consider $RS$-pullback transformations
$RS^2_4\left( 1/3 \atop 3+3+3+3 \right|  {1/2 \atop 2+2+2+2+2+2} \left| 1/2 \atop 7+2+1+1+1 \right)$
of $E(1/3,1/2,0,1/2;t;z)$ and derive solutions of, say, $P_{VI}(1/2,1/2,1/2,-5/2;t)$, $P_{VI}(1/2,1/2,1/2,-1/2;t)$,  $P_{VI}(1/2,1/2,1/2,1/2;t)$. For this, other proper normalization of $\phi_{12}(x)$ has to be used, similarly as other proper normalization of $\varphi_{12}(x)$ was used in \cite[Section 6]{KV3} 
to compute solutions $y_{63}(t_{60})$,  $y_{62}(t_{60})$, $y_{64}(t_{60})$ of the same equations (respectively). Incidentally, pullbacks with respect to $\phi_{12}(x)$ give exactly the same solutions 
$y_{63}(t_{60})$,  $y_{62}(t_{60})$, $y_{64}(t_{60})$ in \cite{KV3} 
of the same three Painlev\'e equations, up to reparametrization $u\mapsto (u+3)/(1-u)$.}
algebraic solutions of, res\-pectively,  $P_{VI}(1/7,1/7,1/7,5/7;t)$, 
$P_{VI}(2/7,2/7,2/7,4/7;t)$ and $P_{VI}(3/7,3/7,3/7,1/7;t)$. 
However, the three solutions are related by Okamoto transformations.
A solution $y_{74}(t_{70})$ of $P_{VI}(1/7,1/7,1/7,5/7;t)$ can be obtained
using  Theorem \ref{kit:method}. Here is a parametrization:
\begin{eqnarray}
t_{70}=\frac{(u-3)^3(u^2+u+2)^2}{2u^3(u^2+7)^2},\qquad 
y_{74}=\frac{(u-1)(u-3)^2(u^2+u+2)}{2u(u^2+3)(u^2+7)}.
\end{eqnarray}
It is related to the parametrization in \cite{HGBAA} via $u\mapsto 1/s$.
Solutions $y_{75}(t_{70})$, $y_{76}(t_{70})$ of, respectively, 
$P_{VI}(2/7,2/7,2/7,4/7;t)$, $P_{VI}(3/7,3/7,3/7,1/7;t)$, can be obtained 
using Theorem \ref{th:llentry}. The same solutions can be obtained
as Okamoto transformations of $y_{74}(t_{70})$.  In the notation of \cite[(2.3)]{KV1}, we have:
\begin{equation}
y_{75}=K_{[-1/7,-1/7,-1/7,\,5/7;\,t_{70}]}\,y_{74}, \qquad 
y_{76}=K_{[1/7,\,1/7,\,1/7,\,5/7;\,t_{70}]}\,y_{74}.
\end{equation}
Here are parametrizations of $y_{75}(t_{70})$ and $y_{76}(t_{70})$:
\begin{eqnarray}
y_{75}\equal-\frac{(u-3)^2(u^2+u+2)^2(u^2+2u+5)}{6u(u+1)(u-1)(u^2+7)},\\ 
y_{76}\equal\frac{(u-1)(u-3)^2(u^2+u+2)(u^4-4u^3-6u^2-28u-11)}{2u(u^2+7)(u^6+21u^4+3u^2+39)}.
\end{eqnarray}
The solution $y_{75}(t_{70})$ is the Kleinian solution of \cite{Bo1}, reparametrized with
$u\mapsto 3s/(s-2)$.

As noticed in \cite{K4},
there are two composite coverings with ramification patterns (\ref{eq:ramif12b}) or (\ref{eq:ramif18}).
They are compositions of Belyi coverings with a quadratic almost Belyi covering:
\begin{eqnarray} \label{eq:comp1}
& R_4(2\;|\;1+1\;|\;1+1)\circ 
R_3(\,\widehat{2}\;|\;2\;|\;1+1\,)\circ R_3(\,\widehat{2}+1\;|\;2+\widehat{1}\;|\;3\,),\\
\label{eq:comp2} & R_4(1+1\;|\;2\;|\;1+1)\circ R_3(3+3+3\;|\;2+2+2+2+\widehat{1}\;|\;7+1+1).
\end{eqnarray}
Here the compositions are from right to left, and the order 2 ramification points of a subsequent 
quadratic covering are indicated by the hat symbol. The algebraic Painlev\'e VI solutions 
are determined by the quadratic almost Belyi coverings. 
The solutions are related (via fractional-linear or Okamoto transformations) to
the solution $y(t)=\sqrt{t}$ of the general equation $P_{VI}(a,b,b,1-a;t)$. 
We specifically have $a=b=1/8$ or $a=b=1/7$ if we apply Theorem \ref{kit:method}
to the two composite coverings. The Belyi coverings are known from 
algebraic transformations of Gauss hypergeometric functions \cite{hhpgtr}. In particular,
an explicit degree 9 covering is given in \cite[(24)]{hhpgtr}.

Beside the indicated coverings, there is exactly one covering (up to fractional-linear transformations)
to pullback hyperbolic hypergeometric equations. It has ramification pattern (\ref{eq:ramif12b}):
\begin{equation}
\psi_{12}(x)=-\frac{4\left(9x^4+18x^3+3(2s+5)x^2-2(s-2)x+s(s-2)\right))^3}
{(4s+1)^3\left(9x^4+14x^3+3(2s+3)x^2-6sx+s^2\right)},
\end{equation}
To get a proper normalization or apply 
Theorem \ref{th:llentry},  we need to choose the point $x=\infty$ appropriately;
hence first a transformation
\begin{equation}
s\mapsto -\frac14v^2(3v^2+8v+6),\qquad\qquad  x\mapsto \frac1x-\frac12v^2.
\end{equation}
For a proper normalization, we still need to factor the remaining degree 3 factor polynomial
in the denominator,  and localize the points $x=0$, $x=1$, $x=t$ properly. This is achieved with 
the substitutions
\begin{eqnarray*}
&&v\mapsto 
\frac{(u^2-2)(u^4-4u^3+8u^2+8u+4)}{6u(u^2-2u+2)(u^2+2u+2)}, \\
&&x\mapsto \frac{36u^2(u^4+4)}{(u^2\!+\!2u\!-\!2)(u^4\!+\!8u^2\!+\!4)}\!
\left(\frac{8iu(u^2\!-\!2u\!-\!2)(u^2+2)^3x}{(u^4\!-\!4u^3\!+\!8u^2\!+\!8u\!+\!4)^3}
+\frac{(u^2\!+\!2i)(u^2\!+\!2iu\!+\!2)}{\left(u^2\!+\!2(i\!-\!1)u\!+\!2i\right)^3}\right).\qquad
\end{eqnarray*}
Theorem \ref{kit:method} eventually gives the following solution $y_{81}(t_{80})$ of 
$P_{VI}(1/8,1/8,1/8,7/8;t_{80})$:
\begin{eqnarray}
\!t_{80}\equal 
\frac{i\left(u+i-1\right)^2\left(u-i+1\right)^2
\left(u^2+2(i\!+\!1)u-2i\right)^3\left(u^2-2(i\!+\!1)u-2i\right)^3}
{64u^2\left(u^2-2\right)^3\left(u^2+2\right)^3},\\
\label{eq:y81}
\!y_{81}\equal \! -i\frac{(u\!+\!i\!-\!1)(u\!-\!i\!+\!1)
(u^2\!+\!2iu\!+\!2)(u^2+2(i\!+\!1)u-2i)^2(u^2-2(i\!+\!1)u-2i)}
{8u\left(u^2-2\right)^2\left(u^2+2\right)\left(u^2-2u-2\right)}.\qquad
\end{eqnarray}
The numerator of $t_{80}$ can also be written as $i\left(u^2+2i\right)^2\left(u^4-12iu^2-4\right)^3$, 
for instance. 
The substitution $u\mapsto -(1+i)/u$ gives the parametrization \cite[(4.12)--(4.13)]{K4}.
 
With the same proper normalization of $\psi_{12}(x)$, one may consider $RS$-transformations
$RS^2_4\left( 1/3 \atop 3+3+3+3 \right|  {1/2 \atop 2+2+2+2+2+2} \left| 3/8 \atop 8+1+1+1+1 \right)$,\,
$RS^2_4\left( 1/3 \atop 3+3+3+3 \right|  {1/2 \atop 2+2+2+2+2+2} \left| 1/4 \atop 8+1+1+1+1 \right)$ 
and $RS^2_4\left( 1/3 \atop 3+3+3+3 \right|  {1/2 \atop 2+2+2+2+2+2} \left| 1/2 \atop 8+1+1+1+1 \right)$
to derive solutions of  $P_{VI}(3/8,3/8,3/8,5/8;t_{80})$,  $P_{VI}(1/4,1/4,1/4,1/4;t_{80})$
and $P_{VI}(1/2,1/2,1/2,1/2;t_{80})$, respectively. The last two equations turn out to be 
the same as $y_{62}(t_{60})$ and $y_{64}(t_{60})$, respectively; Theorem \ref{th:llentry} 
gives expressions reparametrized by 
\begin{equation} \label{eq:tr48}
u\to -\frac{u^4+12iu^2-4}{u^4-4iu^2-4}. 
\end{equation}
The solution of $P_{VI}(3/8,3/8,3/8,5/8;t_{80})$ is
\begin{eqnarray} \label{eq:y84}
y_{83}\equal \frac{i\, (u\!+\!i\!-\!1)(u\!-\!i\!+\!1) 
\left(u^2+2(i\!+\!1)u-2i\right)^2\left(u^2-2(i\!+\!1)u-2i\right)}
{8u\left(u^2-2\right)^2\left(u^2+2\right)\left(u^6+6u^5+6u^4+16u^3-12u^2+24u-8\right)} \nonumber\\
 && \times\left(u^6-6iu^5-6u^4+16iu^3-12u^2-24iu+8\right).
\end{eqnarray}
This solution is presented in \cite[pg.~102]{Bo3}, reparametrized with $u\mapsto(i-1)s$. 
The same solution can be obtained by an Okamoto transformation: 
$y_{83}=K_{[-1/8,-1/8,-1/8,\,7/8;\,t_{80}]}\,y_{81}$.

As was suspected in \cite{K4}, the solutions $y_{62}(t_{60})$ of $P_{VI}(1/4,1/4,1/4,1/4;t_{80})$
and $y_{81}(t_{80})$ of $P_{VI}(1/8,1/8,1/8,7/8;t_{80})$ are related by a sequence of two quadratic transformations. Indeed, a fractional-linear transformation of $K_{[-1/4,-1/4,-1/4,\,1/4;\,t_{60}]}\,y_{62}$
solves $P_{VI}(0,0,1/2,1;t_{80})$, and then we can apply the following result on composition of two quadratic transformations. After substitution (\ref{eq:tr48}) the square roots are extractable;
see Lemma \ref{th:quadr2}.

\section{Appendix}
\label{sec:appendix}

Here we briefly recall or consider the following topics:
\begin{itemize}
\item A formula for composition of two subsequent quadratic transformations of Painlev\'e VI functions;
see Lemma \ref{th:quadr2}. 
\item A general formula for the degree of almost Belyi coverings relevant to algebraic Painlev\'e VI solutions; see Lemma \ref{th:cdegree}.
\item A geometric interpretation of the degree formula.
\end{itemize}


\begin{lemma} \label{th:quadr2}
Suppose that $y(t)$ is a solution of $P_{VI}(0,0,a,1;t)$. Then the following expression 
is a solution of $P_{VI}(a/4,a/4,a/4,1-a/4;t)$:
\begin{equation}
t\;\frac{\sqrt{(y-1)(t-1)}+\sqrt{y\,t}+1}{\sqrt{y\,t}+t}.
\end{equation}
\end{lemma}
\begin{proof} 
The result \cite{RGT} of Ramani-Gramatikos-Tamizhmani states that if $Y_0(T_0)$ is
a solution of $P_{VI}(0,b,c,1;T_0)$, and
\begin{equation}
Y_1 = \frac{(\sqrt{Y_0}+1)(\sqrt{T_0}+1)}{(\sqrt{Y_0}-1)(\sqrt{T_0}-1)}, \qquad
T_1 = \frac{(\sqrt{T_0}+1)^2}{(\sqrt{T_0}-1)^2},
\end{equation}
then $Y_1(T_1)$ is a solution of $P_{VI}(b/2,c/2,c/2,1-b/2;T_1)$. We can transform $y(t)$ to a solution
of $P_{VI}(0,a/2,a/2,1;\ldots)$, and then apply the same transformation to get the asserted solution. 
(Other branches of the transformed solution can be obtained by flipping the sign of the square roots 
$\sqrt{(y-1)(t-1)}$ and $\sqrt{y\,t}$.)
\end{proof}

The following is a degree formula for pullback coverings 
generating algebraic Painlev\'e VI solutions by Theorem \ref{kit:method}. 
In particular, it implies that the pullback covering for an icosahedral \cite{Bo2} solution of
$P_{VI}(\nu_0,\nu_1,\nu_t,\nu_\infty;t)$ with $\nu_0,\nu_1,\nu_t,\nu_\infty\in (0,1)$
has the degree $30(\nu_0+\nu_1+\nu_t-\nu_\infty)$.
\begin{lemma} \label{th:cdegree}
In  the situation of Theorem $\ref{kit:method}$, we have, 
if $\frac1{k_0}+\frac1{k_1}+\frac{1}{k_\infty}\neq 1$:
\begin{equation}
\deg\varphi=\left(\frac{a_0}{k_{\varphi(0)}}+\frac{a_1}{k_{\varphi(1)}}+
\frac{a_t}{k_{\varphi(t)}}+\frac{a_{\infty}}{k_{\varphi(\infty)}}-1\right)
\left/\left(\frac1{k_0}+\frac1{k_1}+\frac{1}{k_\infty}-1\right)\right..
\end{equation}
\end{lemma}
\begin{proof}
Let $d$ denote the degree of $\varphi$. Let $b_0$, $b_1$ respectively $b_{\infty}$ denote the sums of those $a_x$ with $x\in\{0,1,t,\infty\}$ such that, respectively, $\varphi(x)=0$, $\varphi(x)=1$, 
$\varphi(x)=\infty$.  By the Hurwitz formula, we have
\begin{eqnarray*}
2d-2\equal (k_0-1)\frac{d-b_0}{k_0}+(k_1-1)\frac{d-b_1}{k_1}+(k_\infty-1)\frac{d-b_\infty}{k_\infty}\\
&& +(a_0-1)+(a_1-1)+(a_t-1)+(a_\infty-1)+1.
\end{eqnarray*}
The formula follows, since $b_0+b_1+b_\infty=a_0+a_1+a_t+a_\infty$.  
\end{proof}

Notice that this Lemma implies that it is not possible to obtain solutions like $y_{72}(t_{12})$
$y_{75}(t_{70})$ using Theorem \ref{kit:method}: the degree of the covering would be negative.
In other words, we cannot pullback the hyperbolic hypergeometric equation $E(1/3,1/2,0,1/7;t;z)$ 
to  the equations like $E(2/7,2/7,2/7,1/3;y_{72};z)$ or $E(2/7,2/7,2/7,4/7;y_{75};z)$. 
As one can see, there are just a few pullback coverings for infinitely many ``hyperbolic" Painlev\'e VI solutions. This is in contrast to icosahedral  Painlev\'e VI solutions (or more generally, solutions corresponding to  Fuchsian systems with a finite monodromy), which can be obtained from a standard icosahedral hypergeometric equation thanks to Klein's theorem. 

There is a geometric interpretation of this degree formula. 
If $\frac1{k_0}+\frac1{k_1}+\frac{1}{k_\infty}>1$, then
the expression $\left(\frac1{k_0}+\frac1{k_1}+\frac{1}{k_\infty}-1\right)\!\pi$
is the area of the spherical triangle with the angles $\pi/k_0$, $\pi/k_1$, $\pi/k_\infty$
in the standard Riemannian metric on the sphere. The spherical triangle is the image of
the upper-half plane of a Schwarz map for a hypergeometric differential equation with the local exponent differences $1/k_0$, $1/k_1$, $1/k_\infty$. The image of a Schwarz map for a scalar Fuchsian equation associated with (\ref{kit:fuchs}) is a degenerate pentagon, with four angles equal to
$a_0\pi/k_{\varphi(0)}$, $a_1\pi/k_{\varphi(1)}$, $a_t\pi/k_{\varphi(t)}$, 
$\pi-a_\infty\pi/k_{\varphi(\infty)}$, and one angles (corresponding to
the extra ramification point) equal to $2\pi$. The area of the degenerate pentagon is equal to
$\left(\frac{a_0}{k_{\varphi(0)}}+\frac{a_1}{k_{\varphi(1)}}+\frac{a_t}{k_{\varphi(t)}}
+\frac{a_{\infty}}{k_{\varphi(\infty)}}-1\right)\!\pi$. 
If the covering $z=\varphi(x)$ can be defined over $\RR$,
then the degenerate pentagon can be triangulated into the Schwarz triangles with the angles 
$\pi/k_0$, $\pi/k_1$, $\pi/k_\infty$, respecting analytic continuation 
(between the two complex half-planes) in the fiber (with respect to $\varphi$)
of the degenerate pentagon. If $\frac1{k_0}+\frac1{k_1}+\frac{1}{k_\infty}<1$ then we have hyperbolic triangles instead of spherical triangles, with the area 
$\left(1-\frac1{k_0}-\frac1{k_1}-\frac{1}{k_\infty}\right)\!\pi$ with respect to a hyperbolic metric,
but other features are the same.

\begin{figure}
\[
\setlength{\unitlength}{0.75pt}
\begin{picture}(510,388)(15,-80)
\put(19,-20){\resizebox{7.5cm}{!}{\includegraphics{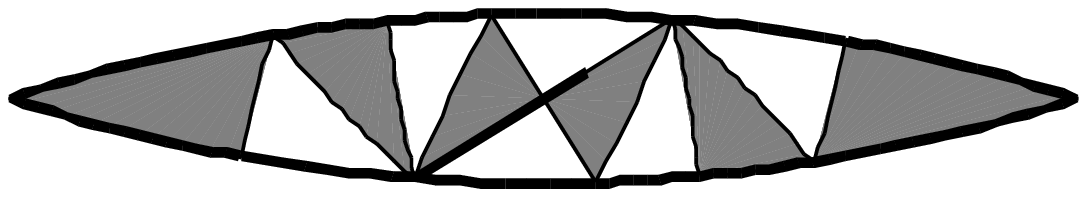}}}
\put(260,-154){\resizebox{7cm}{!}{\includegraphics{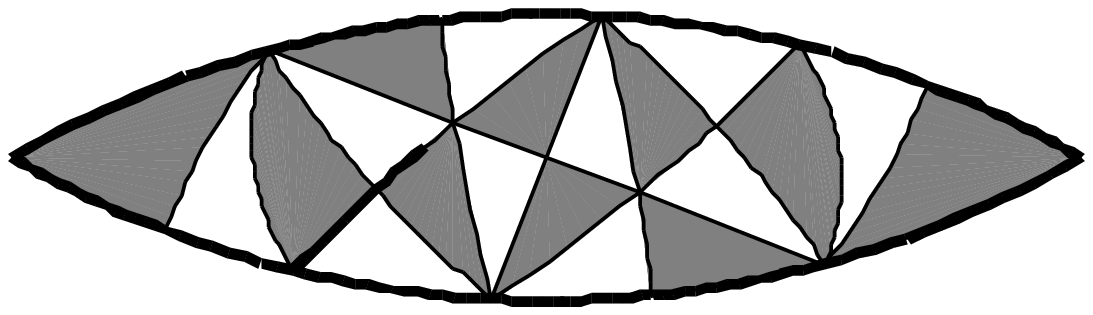}}}
\put(305,8){\resizebox{5.8cm}{!}{\includegraphics{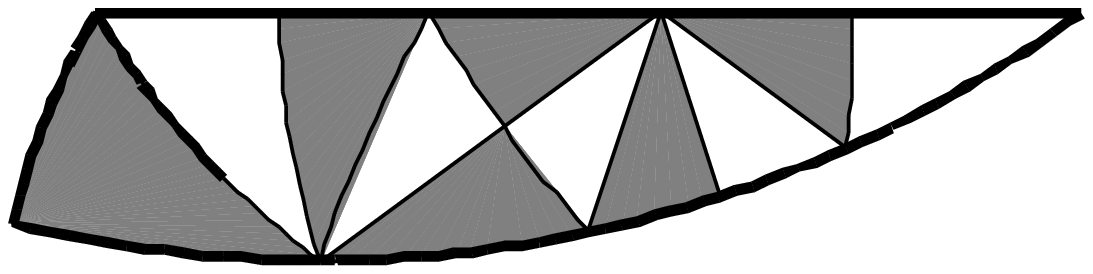}}}
\put(345,148){\resizebox{5.8cm}{!}{\includegraphics{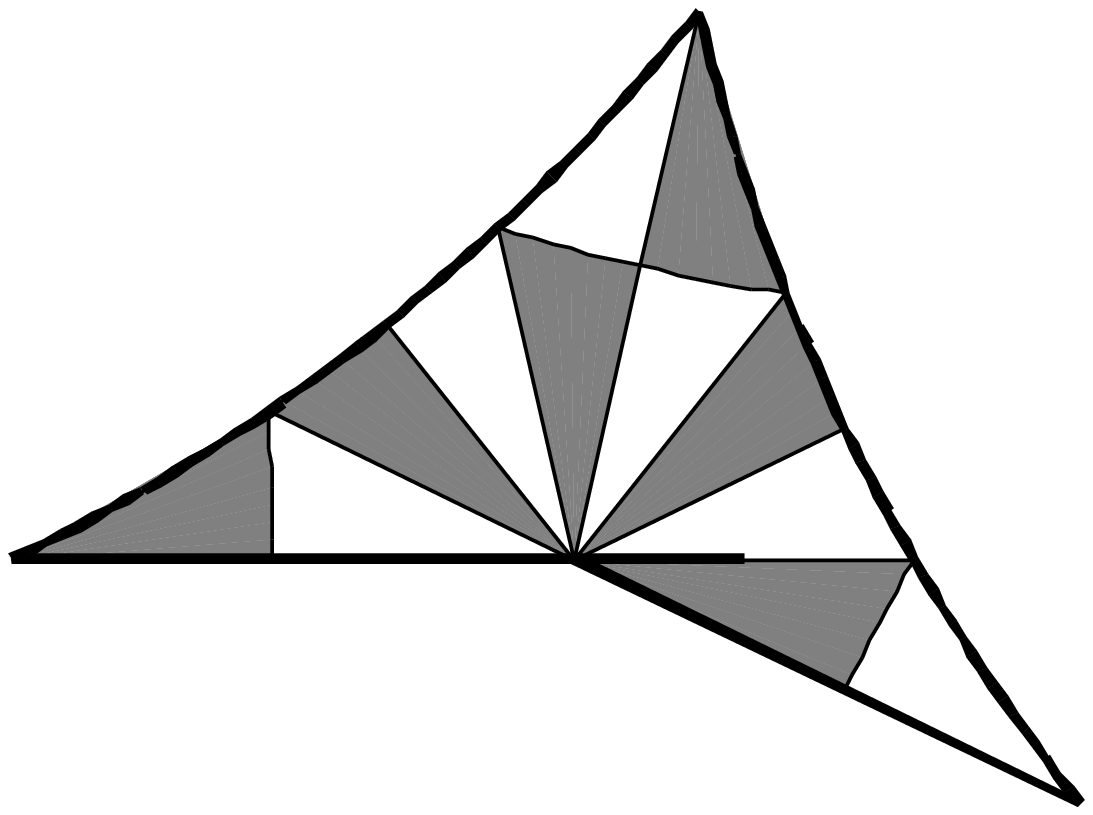}}}
\put(180,140){\resizebox{5cm}{!}{\includegraphics{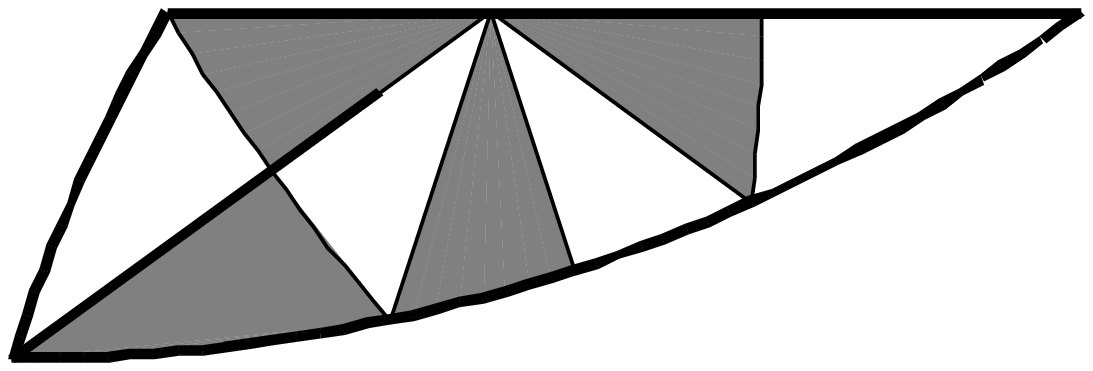}}}
\put(-16,155){\resizebox{5.6cm}{!}{\includegraphics{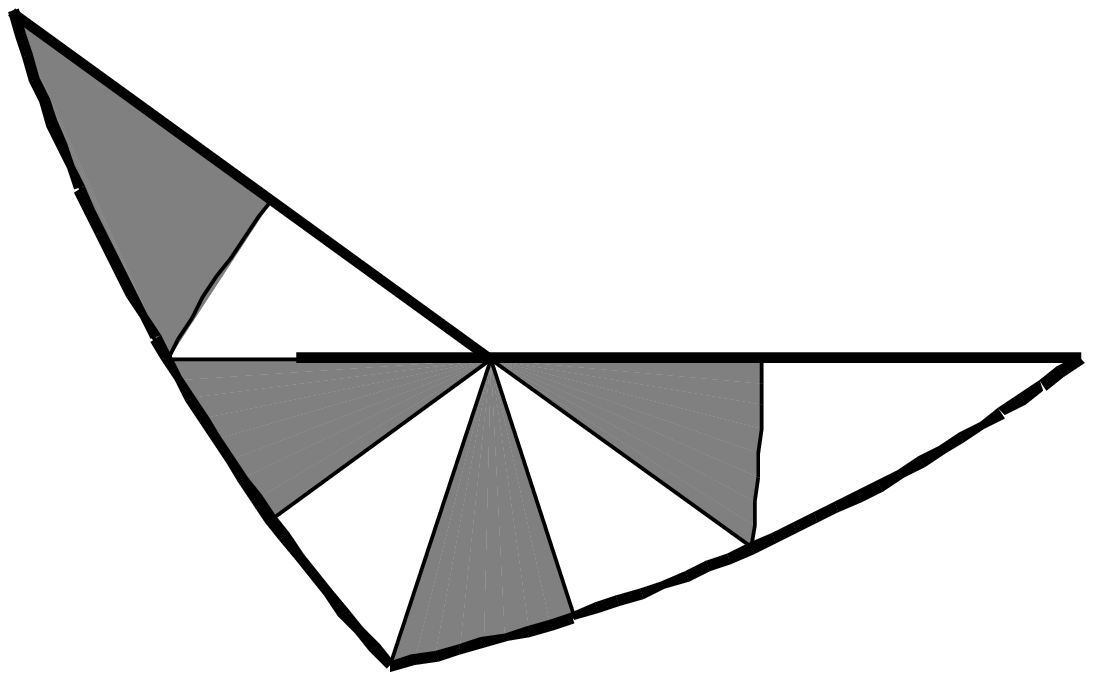}}}
\put(36,-100){\resizebox{4.9cm}{!}{\includegraphics{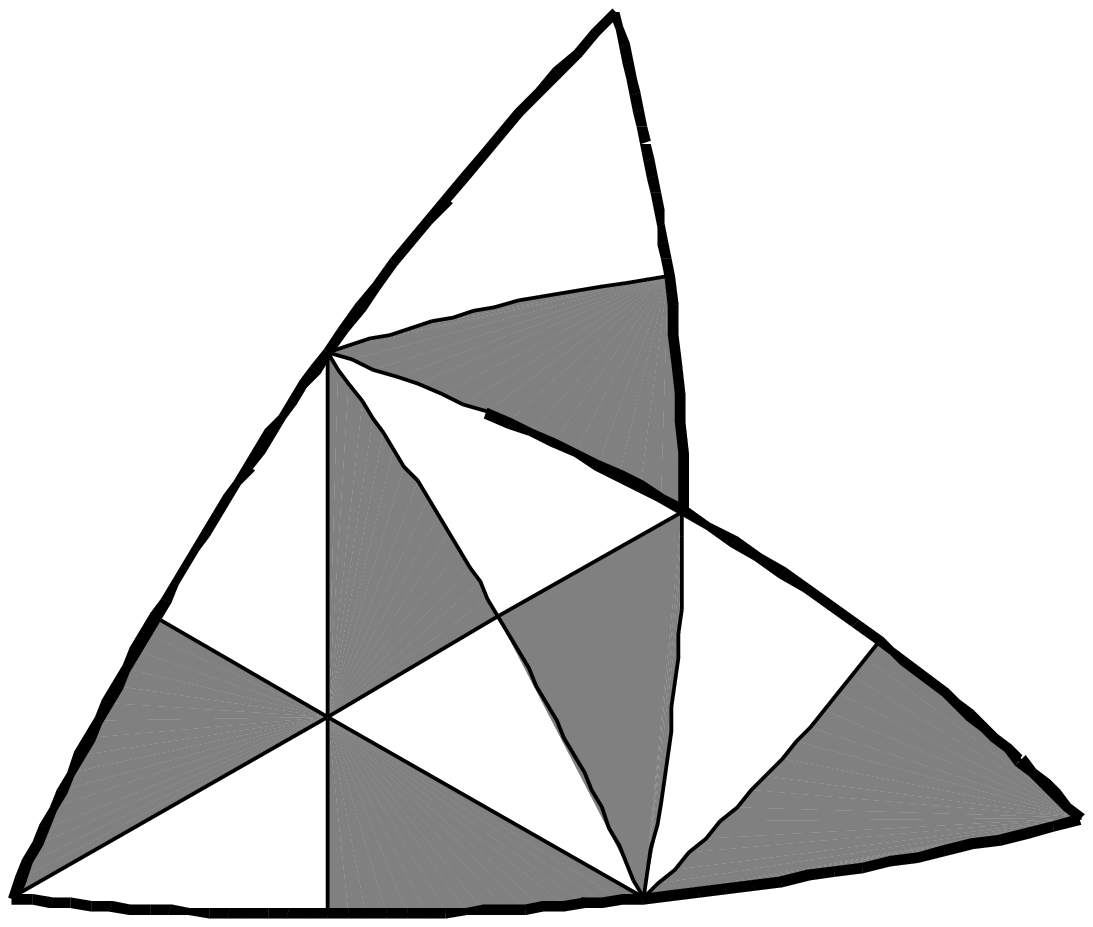}}}
\put(120,195){\it(a)} \put(294,195){\it(b)} \put(460,195){\it(c)}
\put(256,85){\it(d)} \put(460,85){\it(e)} \put(210,-69){\it(f)} \put(460,-69){\it(g)}
\end{picture}
\]
\caption{Triangulations for Schwarz maps}
\label{triangles237}
\end{figure}

Figures \ref{triangles237}{\em (a)} and {\em (b)} depict Schwarz triangulations for the degree
8 map $\widehat{\varphi}_8(x)$ in \cite[(2.7)]{KV3}. The cut for the fifth vertex in Figure
\ref{triangles237}{\em (b)} can either include or do not reach the interior vertex. 
Two different figures correspond to two connected components over
$\RR$ of the Hurwitz curve $w^2=s(s-1)(s+3)(s+8)$. The two components can be distinguished by 
the cut from a point above $z=0$: in Figure \ref{triangles237}{\em (a)} the cut goes towards a point
above  $z=\infty$, while in Figure \ref{triangles237}{\em (b)} it goes towards a point above $z=1$.
One can evaluate $\widehat{\varphi}_8(x)$ at the extra ramification point: 
\begin{equation}
\widehat{\varphi}_8(y_{26})=
-\frac{3125(u+3)(u+2)^4(2u+1)^2u^2(u-1)^3}{4(u+8)(u^3+4u^2+2u+2)^5(u-2)^2}.
\end{equation}
The value $\widehat{\varphi}_8(y_{26})$ oscillates between $z=0$ and $z=1$ for $u\in[-3,0]$, 
and the value is negative or $z=0$, $z=\infty$ when $u\ge 1$ or $u\le -8$. Hence, 
Figure \ref{triangles237}{\em (b)} corresponds to the real component with $u\in[-3,0]$, and 
Figure \ref{triangles237}{\em (a)} corresponds to the other real component.
Notice that  $\widehat{\varphi}_8(y_{26})$, as a function of $u$, is a Belyi map.

Figure \ref{triangles237}{\em (d)} depicts a Schwarz triangulation for the degree 12 map $\widehat{\varphi}_{12}(x)$ in \cite[(2.13)]{KV3}. Figure \ref{triangles237}{\em (c)} depicts 
a Schwarz triangulation for a normalization of ${\phi}_{12}(x)$ here; 
this is a hyperbolic triangulation. Schwarz triangulations for our $\widehat{\varphi}_{10}(x)$ and
normalized composite coverings for (\ref{eq:comp1}) are modifications of two triangulations 
for Belyi coverings  in \cite[Fig. 1]{hhpgtr}: there has to be a cut from the vertices with
the angles $2\pi/7$ and $2\pi/8$.  
Figures \ref{triangles237}{\em (e)}, {\em (f)}, {\em (g)} depict Schwarz triangulations for the degree 11, 12, 20 maps in \cite{KV2}. Note that the lens shaped figures {\em (d)} and {\em (g)} correspond precisely to Dubrovin-Mazzocco solutions. 

Not all almost Belyi coverings have Schwarz tringulations. If a covering is not defined over $\RR$,
analytic continuations of Schwarz maps for the original and transformed equations do not match.
For example, normalizations of ${\psi}_{12}(x)$ can be defined only over $\QQ(i)$. 
Normalized composite coverings for (\ref{eq:comp1}), or the composite degree 20 map 
$\varphi_4\circ\varphi_5(x)$ in \cite[Section 5]{KV2} are not defined over $\RR$ either;
nor they have Schwarz triangulations.

\small

\bibliographystyle{alpha}

\begin{thebibliography}{99}
\bibitem{AK1}
F.~V.~Andreev and A.~V.~Kitaev, Some Examples of $RS^2_3(3)$-Transformations
of Ranks $5$ and $6$ as the Higher Order Transformations for the
Hypergeometric Function, {\it Ramanujan J.} {\bf 7} (2003), no.~4, 455-476,
(http://xyz.lanl.gov, nlin.SI/0012052, 1-20, 2000).
\bibitem{AK2}
F.~V.~Andreev and A.~V.~Kitaev, Transformations ${RS}_4^2(3)$ of the Ranks
$\leq4$ and Algebraic Solutions of the Sixth Painlev\'e Equation, {\it Comm.
Math. Phys.} {\bf 228} (2002), 151--176, (http://xyz.lanl.gov,
nlin.SI/0107074, 1--26, 2001).
\bibitem{BeukersWaall}
F.~Beukers, A.~van der Waall, Lam\'e  equations with algebraic solutions,
{\it J. Differential Equations} {\bf 197} (2004), no.~1, 1--25.
\bibitem{Bo1}
P.~Boalch, From Klein Solution Painlev\'e via Fourier, Laplace and
Jimbo, 
{\em Proc. London Math. Soc. (3)}, {\bf 90} (2005), 167-208.
\bibitem{Bo2}
P.~Boalch, The fifty-two icosahedral solutions to Painlev\'e VI, 
{\em J. Reine Angew. Math.} 596 (2006), 183-214. 
\bibitem{Bo3}
P.~Boalch,
{\it Some explicit solutions to the Riemann-Hilbert problem}, in
``Differential Equations and Quantum Groups", 
IRMA Lectures in Mathematics and Theoretical Physics, Vol. 9 (2006),
pg 85--112.
\bibitem{Do}
Ch.~F.~Doran, Algebraic and Geometric Isomonodromic Deformations, {\it J.
Differential Geometry} {\bf 59} (2001), 33-85.
\bibitem{DM}
B.~Dubrovin and M.~Mazzocco, Monodromy of Certain Painlev\'e--VI
Transcendents and Reflection Groups, {\it Invent. Math.} {\bf 141} (2000),
55--147.
\bibitem{FN}
H.~Flashka, A.~C.~Newell,
Monodromy and spectrum preserving deformations. {\em I Commun. Math. Phys.},
v. 76 pp. 67--116 (1980)
\bibitem{F2}
R.~Fuchs, \"Uber lineare homogene Differentialgleichungen zweiter Ordnung
mit drei im Endlichen gelegenen wesentlich singul\"aren Stellen, {\it Math.
Ann.} {\bf70} (1911), 525-549.
\bibitem{JM}
M.~Jimbo and T.~Miwa, Monodromy preserving deformation of linear ordinary
differential equations with rational coefficients II, {\it Physica} {\bf 2D}
(1981), 407--448.
\bibitem{K2}
A.~V.~Kitaev, Special Functions of the Isomonodromy Type, Rational
Transformations of Spectral Parameter, and Algebraic Solutions of the Sixth
Painlev\'e Equation (Russian), {\it Algebra i Analiz} {\bf 14} (2002), o.~3,
121--139. English Translation in {\it St. Petersburg Math. J.} {\bf14},
no.~3, 453--465 (2003). Available at {\sf http://xxx.lanl.gov,
nlin.SI/0102020}.
\bibitem{HGBAA}
A.~V.~Kitaev, Grothendieck's Dessins d'Enfants, Their Deformations and
Algebraic Solutions of the Sixth Painlev\'e and Gauss Hypergeometric
Equations, {\it Algebra i Analiz} {\bf17}, no.~1 (2005), 224-273.
Available at {\sf http://xxx.lanl.gov, nlin.SI/0309078}.
\bibitem{K4}
A.~V.~Kitaev, {\it Remarks Towards Classification of
$RS_4^2(3)$-Trans\-for\-ma\-tions and Algebraic Solutions of the
Sixth Painlev\'e Equation},
available at {\sf http://xxx.lanl.gov, math.CA/0503082}.
\bibitem{Klein84}
F.~Klein, Vorlesungen \"uber das Ikosaedar, B.~G.~Teubner,Leipzig, 1884.
\bibitem{LosTyk}
O.~Lisovyy, Y.~Tykhyy, {\em Algebraic solutions of the sixth Painlev\'e equation}.
Available at {\sf http://arxiv.org/abs/0809.4873}, (2008).
\bibitem{Litcanu}
R.~Litcanu, Lam\'e operators with finite monodromy - a combinatorial approach,
{\it J. Differential Equations} {\bf 207} (2004), 93--116.
\bibitem{Oh}
Y.~Ohyama, S.~Okumura, {\em R.~Fuchs problem of the Painlev\'e Equations
from the first to the fifth}. Available at {\sf http://www.arxiv.org/math.CA/0512243},
(2005).
\bibitem{O}
K.~Okamoto, Studies on the Painlev\'e Equations. I. Sixth Painlev\'e
Equation $P_{VI}$, {\it Annali Mat. Pura Appl.} {\bf 146} (1987), 337--381.
\bibitem{RGT}
A.~Ramani, B.~Grammaticos, T. Tamizhmani, Quadratic relations
in continuous and discrete Painlev\'e Equations. {\it J. Phys. A.: Math.
Gen.} {\bf 33} (2000), 3033--3044.
\bibitem{KV1}
R.~Vidunas, A.~V.~Kitaev, Quadratic Transformations of the Sixth Painlev\'e
Equation. Accepted by {\em Mathematische Nachrichten}. Available
at {\sf http://arxiv.org/abs/math.CA/0511149}.
\bibitem{KV2}
R.~Vidunas, A.~V.~Kitaev, Computation of highly ramified coverings, Available
at {\sf http://arxiv.org/abs/0705.3134}.
\bibitem{KV3}
R.~Vidunas, A.~V.~Kitaev, Schlesinger transformations
for algebraic Painlev\'e VI solutions. Available
at {\sf http://arxiv.org/abs/0810.2766}.
\bibitem{V}
R.~Vidunas, Algebraic Transformations of Gauss Hypergeometric
Functions, submitted to "Funkcialaj Ekvacioj"; preprint
{\sf http://www.arxiv.org/math.CA/0408269} (2004). 
\bibitem{hhpgtr}
R.~Vidunas, Transformations of some Gauss hypergeometric functions, {\it
J. Comp. Appl. Math.}~{\bf178}~(2005), 473--487.
\bibitem{zvonkin}
A.~Zvonkin, {\em Megamaps: Construction and Examples}, Discrete Mathematics
and Theoretical Computer Science Proceedings AA (DM-CCG),
2001, pg. 329-340.
\end{thebibliography}

\end{document}